\documentclass[12pt,reqno]{amsart}
\usepackage{amssymb}
\usepackage{eurosym}
\usepackage{epsfig}
\usepackage{amsmath}
\usepackage{amsfonts}
\usepackage{pgf,tikz}
\usepackage{enumerate}
\usepackage{cite}
\usepackage[colorlinks=true, linkcolor=black, citecolor=black, urlcolor=black]{hyperref}

\usepackage{graphicx}
\usepackage{caption}
\usepackage{subcaption}

\newcommand{\R}{{{\Bbb R}}}

\newtheorem{theorem}{\sc Theorem}[section]

\newtheorem{remark}{\sc Remark}[section]
\newtheorem{corollary}{\sc Corollary}[section]
\newtheorem{example}{\sc Example}[section]

\def\qed{\hbox to 0pt{}\hfill$\rlap{$\sqcap$}\sqcup$\medbreak}

\title[A hybrid Krasnosel'ski\u{i}-Schauder FPT]{A hybrid Krasnosel'ski\u{i}-Schauder fixed point theorem for systems}

\begin{document}
\author[G. Infante]{Gennaro Infante}
\address{Gennaro Infante, Dipartimento di Matematica e Informatica, Universit\`{a} della
Calabria, 87036 Arcavacata di Rende, Cosenza, Italy}%
\email{gennaro.infante@unical.it}%

\author[G. Mascali]{Giovanni Mascali}
\address{Giovanni Mascali, Dipartimento di Matematica e Informatica, Universit\`{a} della
Calabria, 87036 Arcavacata di Rende, Cosenza, Italy}%
\email{giovanni.mascali@unical.it}%

\author[J. Rodr\'iguez--L\'opez]{Jorge Rodr\'iguez--L\'opez}
\address{Jorge Rodr\'iguez--L\'opez, CITMAga \& Departamento de Estat\'{\i}stica, An\'alise Matem\'atica e Optimizaci\'on, Universidade de Santiago de Compostela,  15782, Facultade de Matem\'aticas, Campus Vida, Santiago, Spain}
\email{jorgerodriguez.lopez@usc.es}%

\date{}

\begin{abstract}
We provide new results regarding the localization of the solutions of nonlinear operator systems. We make use of a combination of Krasnosel'ski\u{\i} cone compression-expansion type methodologies and Schauder-type ones. In particular we establish a localization of the solution of the system within the product of a conical shell and of a closed convex set. By iterating this procedure we prove the existence of multiple solutions. We illustrate our theoretical results by applying them to the solvability of systems of Hammerstein integral equations. In the case of two specific boundary value problems and with given nonlinearities, we are also able to obtain a numerical solution, consistent with our theoretical results.  
\end{abstract}


 \maketitle

\section{Introduction}

The Krasnosel'ski\u{\i} compression-expansion fixed point theorem in cones and the Schauder fixed point theorem in normed spaces are among the most well-known and widely employed topological fixed point theorems in the literature. Both results have been extensively applied in order to obtain existence and localization of solutions for a huge variety of boundary value problems and integral equations.    

Our purpose is to present a novel fixed point theorem for operator systems of the form 
\begin{equation}\label{syst}
\left\{\begin{array}{l} u_1=T_1(u_1,u_2), \\ u_2=T_2(u_1,u_2), \end{array} \right. 
\end{equation}
which combines, in a component-wise manner, the assumptions of the Krasnosel'ski\u{\i} cone fixed point theorem with those in the Schauder theorem. Roughly speaking, we assume that for each fixed $u_2$, the operator $T_1(\cdot,u_2)$ satisfies compression-expansion type conditions in the line of Krasnosel'ski\u{\i}'s result and, moreover, for each fixed $u_1$, the operator $T_2(u_1,\cdot)$ fulfills the hypotheses of the Schauder fixed point theorem (see Theorem~\ref{th_KrasSchauder} below). As a result, we obtain a solution $(\bar{u}_1,\bar{u}_2)$ of the system \eqref{syst} with the first component, $\bar{u}_1$, localized in a conical shell and its second one, $\bar{u}_2$, in a certain closed convex set. By an iteration of this methodology we show that it is also possible to prove the existence of multiple solutions.
We emphasize that the localization of solutions of differential systems (or more in general of nonlinear operator equations) plays a key role in various settings, for example when modelling biological or medical phenomena~\cite{Murray}.

The proof relies on an application of the classical Leray-Schauder fixed point index, see \cite{amann,GraDug, guolak}. 
 Our abstract theory {follows the path of previous results for operator systems due to Precup and co-authors \cite{Benedetti,PrecupFPT,PrecupSDC}, who were, in their turn, motivated by the papers~\cite{Avramescu, Perov}. In particular, we benefit from the setting in \cite{PrecupFPT,PrecupSDC}}, where the author established the \textit{vector version of Krasnosel'ski\u{\i} fixed point theorem}, so-called because compression-expansion type conditions are imposed independently in each component of the system~\eqref{syst}. An alternative approach to this fixed point theorem based on fixed point index techniques has been recently presented in \cite{JRL} and complementary results in this line can also be found in \cite{INOP,ima,imap,CPR,Lan1,PreRod}. 

Finally, to show that this hybrid fixed point theorem is a useful tool for the study of systems of nonlinear differential and integral equations, we apply it to systems of Hammerstein integral equations. As a consequence, we obtain existence and localization results for second-order systems subject to Sturm-Liouville type boundary conditions. We illustrate, in the case of two specific examples, the constants that occur in our theoretical results and we also exhibit numerical solutions that are consistent with our theory.
Our results are new and complement previous results in the literature; in particular the recent results of~\cite{BGK23}, where the approach relies on a combination of monotone and variational techniques.

\section{Fixed point theorem for operator systems}\label{sec2}

In the sequel, we need the following notions. A closed convex subset $K$ of a normed linear space $(X,\left\|\cdot\right\|)$ is a \textit{cone} if $\lambda\,u\in K$ for every $u\in K$ and for all $\lambda\geq 0$, and $K\cap (-K)=\{0\}$.

A cone $K$ induces the partial order in $X$ given by $u\preceq v$ if and only if $v-u\in K$. Moreover, we shall say that $u\prec v$ if $v-u\in K\setminus\{0\}$. The cone $K$ is called \textit{normal} if there exists $d>0$ such that $\left\|u\right\|\leq d\left\|v\right\|$ for all $u,v\in X$ with $0\preceq u\preceq v$.

The following notations will be useful: for given $r,R\in \mathbb{R}_+:=[0,\infty)$, $0<r<R$, we define
\[K_{r,R}:=\{u\in K:r<\left\|u\right\|<R \} \quad \text{ and } \quad \overline{K}_{r,R}:=\{u\in K:r\leq\left\|u\right\|\leq R \}. \]

Now, we recall Krasnosel'ski\u{\i} fixed point theorem \cite{kras} (see also \cite{amann, guolak}). 

\begin{theorem}[Krasnosel'ski\u{\i}]
	Let $(X,\left\|\cdot\right\|)$ be a Banach space, $K$ a cone in $X$ and $r,R\in\R_+$, $0<r<R$.
	
	Consider a compact map $T:\overline{K}_{r,R}\rightarrow K$ (i.e., $T$ is continuous and $T(\overline{K}_{r,R} )$ is relatively compact) satisfying one of the following conditions:
	\begin{enumerate}
		\item[$(a)$] $T(u)\nprec u$ if $\left\|u\right\|=r$ and $T(u)\nsucc u$ if $\left\|u\right\|=R$;
		\item[$(b)$] $T(u)\nsucc u$ if $\left\|u\right\|=r$ and $T(u)\nprec u$ if $\left\|u\right\|=R$.
	\end{enumerate}
	Then $T$ has at least a fixed point $u\in K$ with $r\leq\left\|u\right\|\leq R$.
\end{theorem}

It is commonly said that the operator $T$ is \textit{compressive} if condition $(a)$ in Krasnosel'ski\u{\i} fixed point theorem holds and, \textit{expansive} in case it satisfies $(b)$.	
It is also well-known that conditions $(a)$ and $(b)$ can be weakened as \textit{homotopy type conditions} and the result can be stated for general open subsets of the cone (not necessarily balls). In this way, we have the \textit{homotopy version of Krasnosel'ski\u{\i} theorem} or \textit{Krasnosel'ski\u{\i}-Benjamin theorem}, see for instance~\cite{amann, guolak} and cf.~\cite{GraDug}.

\begin{theorem}[Krasnosel'ski\u{\i}-Benjamin]
	Let $(X,\left\|\cdot\right\|)$ be a {Banach} space, $K$ a cone in $X$ and $U$ and $V$ bounded and relatively open subsets of $K$ such that $0\in V\subset\overline{V}\subset U$.
	
	Assume that $T:\overline{U}\setminus V\rightarrow K$ is a compact map and there exists $h\in K\setminus\{0\}$ such that one of the following conditions is satisfied:
	\begin{enumerate}
		\item[$(a)$] $T(u)+\mu\,h\neq u$ if $u\in\partial\,V$ and $\mu> 0$, and $T(u)\neq \lambda\, u$ if $u\in\partial\,U$ and $\lambda> 1$;
		\item[$(b)$] $T(u)\neq \lambda\, u$ if $u\in\partial\,V$ and $\lambda> 1$, and $T(u)+\mu\,h\neq u$ if $u\in\partial\,U$ and $\mu> 0$. 
	\end{enumerate}
	Then $T$ has at least a fixed point $u\in \overline{U}\setminus V$.
\end{theorem}

For the sake of completeness, let us recall also the Schauder fixed point theorem, see for example~\cite{GraDug}.

\begin{theorem}[Schauder]
	Let $D$ be a non-empty, closed and convex subset of a {Banach} space $(X,\left\|\cdot\right\|)$ and $T:D\rightarrow D$ be a compact map. Then $T$ has at least one fixed point in~$D$.
\end{theorem}

Let $(X,\left\|\cdot\right\|_X)$ and $(Y,\left\|\cdot\right\|_Y)$ be {Banach} spaces and consider the cartesian product $X\times Y$. When no confusion may occur, both norms $\left\|\cdot\right\|_X$ and $\left\|\cdot\right\|_Y$ will be simply denoted by $\left\|\cdot\right\|$.

Our aim is to present a fixed point theorem for systems with Krasnosel'ski\u{\i}-type conditions in one component and Schauder-type conditions in the other one.

\begin{theorem}\label{th_KrasSchauder}
	Let $U$ and $V$ be bounded and relatively open subsets of a cone $K_1$ of the {Banach} space $X$ such that $0\in V\subset\overline{V}\subset U$ and $D$ be a closed convex subset of the {Banach} space $Y$. 
	
	Assume that $T=(T_1,T_2):\left(\overline{U}\setminus V\right)\times D\rightarrow K_1\times D$ is a compact map and there exists $h\in K_1\setminus\{0\}$ such that either of the following conditions holds in $\left(\overline{U}\setminus V\right)\times D$:
	\begin{enumerate}
		\item[$(a)$] $T_1(u)+\mu\,h\neq u_1$ if $u_1\in \partial\,V$ and $\mu> 0$, and $T_1(u)\neq \lambda\, u_1$ if $u_1\in\partial\,U$ and $\lambda>1$; or
		\item[$(b)$] $T_1(u)\neq \lambda\, u_1$ if $u_1\in \partial\,V$ and $\lambda> 1$, and $T_1(u)+\mu\,h\neq u_1$ if $u_1\in\partial\,U$ and $\mu> 0$. 
	\end{enumerate}
	Then $T$ has at least a fixed point $u=(u_1,u_2)\in K_1\times D$ with $u_1\in \overline{U}\setminus V$.
\end{theorem}

\noindent
{\bf Proof.} 
First of all, by the Dugundji extension Theorem (see \cite[Theorem 4.1]{dug}), there exists a continuous map $N=(N_1,N_2):K_1\times D\rightarrow K_1\times D$ such that $N(u)=T(u)$ for all $u\in \left(\overline{U}\setminus V\right)\times D$ and, moreover, $N\left(K_1\times D \right)\subset {\rm co}\,T\left( \left(\overline{U}\setminus V\right)\times D \right)$, where ${\rm co}\,A$ denotes the \textit{convex hull} of $A$, i.e., the smallest convex set containing $A$. According to Mazur's theorem, the set $N\left(K_1\times D \right)$ is relatively compact since it is a subset of the closed convex hull of the relatively compact set $T\left( \left(\overline{U}\setminus V\right)\times D \right)$.

Since the set $C:=K_1\times D$ is a closed convex subset of the {Banach} space $X\times Y$ and $N$ is a compact map, the fixed point index of $N$ with respect to $C$ is well-defined for any (relative) open subset $\mathcal{O}\subset C$ such that $N$ is fixed point free on the boundary of $\mathcal{O}$, see \cite{GraDug}. It will be denoted by $i_C(N,\mathcal{O})$. {Note that the boundary or closure of such set $\mathcal{O}$ is always considered with respect to the relative topology of the set $C$.}
Let us suppose that the operator $T$ is fixed point free on the boundary of the set $\left(\overline{U}\setminus V\right)\times D$ (otherwise the proof is finished), which implies that so is the operator $N$.

Assume now that $T$ satisfies condition $(a)$. Let us consider the relative open sets \[\mathcal{U}:=U\times D \quad \text{ and }\quad \mathcal{V}:=V\times D.\]We shall prove that $i_C(N,\mathcal{U})=1$ and $i_C(N,\mathcal{V})=0$. First, observe that, since the $\partial(D)=\emptyset$ in the relative topology, $\partial\,\mathcal{V}=\{(u_1,u_2)\in C\,:\,u_1\in\partial\,V\}$ and thus condition $(a)$ implies that
\[N_1(u)+\mu\,h\neq u_1 \quad \text{if } u\in\partial\,\mathcal{V} \text{ and } \mu> 0. \]
Hence, $i_C(N,\mathcal{V})=0$. { Indeed, let us assume, on the contrary, that $i_C(N,\mathcal{V})\neq 0$. First, since $V$ is bounded, $N_1(\mathcal{V})$ is relatively compact and $h\neq 0$, there exists $\mu_0>0$ large enough such that $u_1\neq N_1(u)+\mu_0\,h$ for all $u=(u_1,u_2)\in \mathcal{V}$. Now, consider the homotopy $\mathcal{H}:[0,1]\times \overline{\mathcal{V}}\rightarrow C$ given by
	\[\mathcal{H}(t,u)=\left(N_1(u)+\mu_0\, t\, h, N_2(u)\right). \]
By the homotopy property of the index,
\[i_{C}(N,\mathcal{V})=i_{C}(\mathcal{H}(0,\cdot),\mathcal{V})=i_{C}(\mathcal{H}(1,\cdot),\mathcal{V}) \]
which, combined with $i_C(N,\mathcal{V})\neq 0$, implies that there exists $u=(u_1,u_2)\in \mathcal{V}$ such that $u_1= N_1(u)+\mu_0\,h$, a contradiction.	 
}

On the other hand, let us fix $\omega\in D$ and consider the homotopy $H:[0,1]\times\overline{\mathcal{U}}\rightarrow C$ defined as
\[ H(t,u)=\left(t\,N_1(u),t\,N_2(u)+(1-t)\,\omega \right).\]
Note that $u\neq H(t,u)$ for all $u\in\partial\,\mathcal{U}$ and all $t\in[0,1]$. {Indeed, if this is not the case, there exist $u=(u_1,u_2)\in C$ with $u_1\in\partial\,U$ and $t\in[0,1]$ such that $u_1=t\,N_1(u)$ and $u_2=t\,N_2(u)+(1-t)\,\omega$. On the one hand, if $t=0$, then $u_1=0$, which contradicts that $u_1\in\partial\,U$ since $U$ is open and $0\in U$. On the other hand, the operator $T$ is fixed point free on the boundary of the set $\partial (U \times D)\subset\partial \left((\overline{U}\setminus V)\times D\right)$, which implies that $t\neq 1$ (note that the inclusion $\partial (U \times D)\subset\partial \left((\overline{U}\setminus V)\times D\right) $ follows from $\overline{V}\subset U$ and the fact that we are employing the relative topology with respect to $C$). Then there exist $u=(u_1,u_2)\in C$ with $u_1\in\partial\,U$ and $t\in(0,1)$ such that $u_1=t\,N_1(u)$, that is, $(1/t)u_1=T_1(u)$, a contradiction with condition $(a)$.} 
Since $\omega\in D$ and the map $H(0,\cdot)\equiv\omega$, we obtain \[i_C(H(0,\cdot),\mathcal{U})=1.\]

Then, by the homotopy property of the fixed point index, we deduce that
\[ i_C(N,\mathcal{U})=i_C(H(1,\cdot),\mathcal{U})=i_C(H(0,\cdot),\mathcal{U})=1,\]
since $(0,\omega)\in\mathcal{U}$.
Finally, by the additivity property of the index,
\[ i_C(N,\mathcal{U}\setminus\overline{\mathcal{V}})=i_C(N,\mathcal{U})-i_C(N,\mathcal{V})=1,\]
and, by the definition of $N$ we have that $N=T$ on $\mathcal{U}\setminus\overline{\mathcal{V}}$, so we deduce that $T$ has at least one fixed point in $\mathcal{U}\setminus\overline{\mathcal{V}}$. 

The reasoning is analogous if condition $(b)$ is satisfied. In that case we obtain $i_C(N,\mathcal{U})=0$ and $i_C(N,\mathcal{V})=1$. Therefore, it follows from the additivity property of the index that $i_C(N,\mathcal{U}\setminus\overline{\mathcal{V}})=-1$ and thus $T$ has at least one fixed point in $\mathcal{U}\setminus\overline{\mathcal{V}}$.
\qed

\begin{remark}\label{rmk_computation_index}
	From the proof of Theorem \ref{th_KrasSchauder}, one has that if $T$ has no fixed points on $\partial\,(U\setminus\overline{V})\times D$, then $i_{C}(T,(U\setminus\overline{V})\times D )=1$ provided that $T_1$ is a compressive operator (i.e., if $T$ satisfies condition $(a)$) and that $i_{C}(T,(U\setminus\overline{V})\times D )=-1$ if $T$ is expansive (i.e., if condition $(b)$ is fulfilled). These computations of the fixed point index are useful in order to obtain multiple fixed points of the operator $T$.
\end{remark}

{
\begin{remark}
	In view of the proof of Theorem \ref{th_KrasSchauder}, it is not needed that the open set $U$ be bounded in order to prove that $i_{C}(N,\mathcal{U})=1$ (exactly as in the classical case, cf. \cite[Ch. 12, Theorem 7.3]{GraDug}). Nevertheless, to prove that $i_{C}(N,\mathcal{V})=0$, it is employed that $V$ is bounded (again as in the classical case, cf. \cite[Ch. 12, Theorem 7.11]{GraDug}). 
	Therefore, for a compressive operator $T_1$, the result remains valid even if $U$ is unbounded provided that $V$ be bounded.      
\end{remark}
}

Clearly, in the particular case in which $U$ and $V$ are the intersection of two open balls with the cone $K_1$, we obtain the following result.

\begin{corollary}
	Let $r,R\in \mathbb{R}_{+}$ be positive numbers with $r<R$ and $D$ be a closed convex subset of the {Banach} space $Y$. 
	
	Assume that $T=(T_1,T_2):(\overline{K}_{1})_{r,R}\times D\rightarrow K_1\times D$ is a compact map and there exists $h\in K_1\setminus\{0\}$ such that either of the following conditions holds in $(\overline{K}_{1})_{r,R}\times D$:
	\begin{enumerate}
		\item[$(a)$] $T_1(u)+\mu\,h\neq u_1$ if $\left\|u_1\right\|=r$ and $\mu> 0$, and $T_1(u)\neq \lambda\, u_1$ if $\left\|u_1\right\|=R$ and $\lambda> 1$; or
		\item[$(b)$] $T_1(u)\neq \lambda\, u_1$ if $\left\|u_1\right\|=r$ and $\lambda> 1$, and $T_1(u)+\mu\,h\neq u_1$ if $\left\|u_1\right\|=R$ and $\mu> 0$. 
	\end{enumerate}
	Then $T$ has at least a fixed point $u=(u_1,u_2)\in K_1\times D$ with $r\leq\left\|u_1\right\|\leq R$.
\end{corollary}

As a consequence of Theorem \ref{th_KrasSchauder}, we obtain the following criteria when the closed convex set $D$ is an ordered interval $[\alpha,\beta]$ with the order, $\preceq$, given by a cone $K_2$ in the {Banach} space $Y$ {-- that is}, \[[\alpha,\beta]:=\{y\in Y :\alpha\preceq y\preceq \beta \}.\]

\begin{corollary}\label{hybrid_FPT}
	Let $U$ and $V$ be bounded and relatively open subsets of a cone $K_1$ of the {Banach} space $X$ such that $0\in V\subset\overline{V}\subset U$ and $\alpha,\beta\in Y$ such that $\alpha\preceq\beta$. 
	
		Assume that $T=(T_1,T_2):\left(\overline{U}\setminus V\right)\times[\alpha,\beta]\rightarrow K_1\times[\alpha,\beta]$ is a compact map and there exists $h\in K_1\setminus\{0\}$ such that either of the following conditions holds in $\left(\overline{U}\setminus V\right)\times[\alpha,\beta]$:
		\begin{enumerate}
			\item[$(a)$] $T_1(u)+\mu\,h\neq u_1$ if $u_1\in \partial\,V$ and $\mu> 0$, and $T_1(u)\neq \lambda\, u_1$ if $u_1\in\partial\,U$ and $\lambda> 1$; or
			\item[$(b)$] $T_1(u)\neq \lambda\, u_1$ if $u_1\in \partial\,V$ and $\lambda> 1$, and $T_1(u)+\mu\,h\neq u_1$ if $u_1\in\partial\,U$ and $\mu> 0$. 
		\end{enumerate}
	Then $T$ has at least a fixed point $u=(u_1,u_2)\in K_1\times [\alpha,\beta]$ with $u_1\in \overline{U}\setminus V$.
\end{corollary}

\begin{remark}
	If the cone $K_2$ is normal, then the interval $[\alpha,\beta]$ is bounded in $Y$ and thus $\left(\overline{U}\setminus V\right)\times[\alpha,\beta]$ is a bounded subset of $X\times Y$. Hence, the conclusion of Corollary \ref{hybrid_FPT} remains valid if we assume that $T=(T_1,T_2):K_1\times[\alpha,\beta]\rightarrow K_1\times[\alpha,\beta]$ is a completely continuous map (i.e., $T$ is continuous and maps bounded sets into relatively compact ones) satisfying either condition $(a)$ or $(b)$.
\end{remark}

If the operator $T_2$ is increasing with respect to the second variable, the assumption \[T_2\left((\overline{K}_{1})_{r,R}\times[\alpha,\beta] \right)\subset [\alpha,\beta]\] can be clearly deduced from the following one: for all $u_1\in (\overline{K}_{1})_{r,R}$, the order relations $\alpha\preceq T_2(u_1,\alpha)$ and $T_2(u_1,\beta)\preceq \beta$ hold. 

\begin{corollary}
	Let $r,R\in \mathbb{R}_{+}$ be positive numbers with $r<R$ and $\alpha,\beta\in Y$ such that $\alpha\preceq\beta$. 
	
	Assume that $T=(T_1,T_2):(\overline{K}_{1})_{r,R}\times[\alpha,\beta]\rightarrow K_1\times Y$ is a compact map and satisfies the following conditions:
	\begin{enumerate}
		\item One of the following two conditions holds in $(\overline{K}_{1})_{r,R}\times[\alpha,\beta]$:
		\begin{enumerate}
			\item $T_1(u)\nprec u_1$ if $\left\|u_1\right\|=r$  and $T_1(u)\nsucc u_1$ if $\left\|u_1\right\|=R$; or
			\item $T_1(u)\nsucc u_1$ if $\left\|u_1\right\|=r$ and $T_1(u)\nprec u_1$ if $\left\|u_1\right\|=R$. 
		\end{enumerate}
		\item for all $u_1\in (\overline{K}_{1})_{r,R}$, the map $T_2(u_1,\cdot)$ is increasing {and, furthermore,} $\alpha\preceq T_2(u_1,\alpha)$ and $T_2(u_1,\beta)\preceq~\beta$.
	\end{enumerate}
	Then $T$ has at least a fixed point $u=(u_1,u_2)\in K_1\times Y$ with $r\leq\left\|u_1\right\|\leq R$ and $\alpha\preceq u_2\preceq\beta$.
\end{corollary}

Obviously, two distinct fixed points of the operator $T$ can be obtained if the hypotheses of Theorem \ref{th_KrasSchauder} hold in two disjoint domains of the form $\left(\overline{U}_1\setminus V_1\right)\times D_1$ and $\left(\overline{U}_2\setminus V_2\right)\times D_2$. In case of operator systems (of two equations), where by a fixed point we mean a pair $(u_1, u_2)$, two fixed points are different if they differ at least on one of their components, not necessarily on both. Hence, we can look for different fixed points in two disjoint regions of the form $\left(\overline{U}_j\setminus V_j\right)\times D$, $j=1,2$, in such a way that the second component of both fixed points will be localized in the same set $D$.

Furthermore, if we can guarantee that these fixed points are not located on the boundary of the sets $\left(\overline{U}_j\setminus V_j\right)\times D$, then the computation of the fixed point index ensures the existence of a third fixed point.

\begin{theorem}
	Let $U_1$, $U_2$, $V_1$ and $V_2$ be bounded and relatively open subsets of a cone $K_1$ of the {Banach} space $X$ such that $0\in V_1$, $\overline{V}_j\subset U_j$ ($j=1,2$), $\overline{U}_1\subset V_2$ and $D$ be a closed convex subset of the {Banach} space $Y$. 
	
	Assume that $T=(T_1,T_2):\left(\overline{U}_2\setminus V_1\right)\times D\rightarrow K_1\times D$ is a compact map and for each $j\in\{1,2\}$ there exist $h_j\in K_1\setminus\{0\}$ such that the following conditions hold:
	\begin{enumerate}
		\item[$(i)$] $T_1(u)+\mu\,h_j\neq u_1$ if $u_1\in \partial\,V_j$, $u_2\in D$ and $\mu\geq 0$;
		\item[$(ii)$] $T_1(u)\neq \lambda\, u_1$ if $u_1\in\partial\,U_j$, $u_2\in D$ and $\lambda\geq 1$.
	\end{enumerate}
	Then $T$ has at least three distinct fixed points $(u_1,u_2),(v_1,v_2),(w_1,w_2)\in K_1\times D$ with $u_1\in U_1\setminus \overline{V}_1$, $v_1\in U_2\setminus \overline{V}_2$ and $w_1\in V_2\setminus \overline{U}_1$.
\end{theorem}

\noindent
{\bf Proof.} By the computations of the fixed point index, we have
\[i_{C}(T,(U_1\setminus\overline{V}_1)\times D )=i_{C}(T,(U_2\setminus\overline{V}_2)\times D )=1, \quad i_{C}(T,(V_2\setminus\overline{U}_1)\times D )=-1. \]
Therefore, the conclusion follows from the existence property of the fixed point index.
\qed

\section{Application to Hammerstein systems}

Consider the following system of Hammerstein type equations
\begin{equation}\label{eq_Ham}
	\begin{array}{r}
		u(t)=\displaystyle\int_{0}^{1} k_1(t,s)f(s,u(s),v(s))\,ds:=T_1(u,v)(t), \\[0.3cm]
		v(t)=\displaystyle\int_{0}^{1} k_2(t,s)g(s,u(s),v(s))\,ds:=T_2(u,v)(t),
	\end{array}
\end{equation}
where $I:=[0,1]$ and the following assumptions are satisfied:
\begin{enumerate}
	\item[$(H_1)$] the kernels $k_1,k_2:I^2\rightarrow\mathbb{R}_+$ are continuous;
	\item[$(H_2)$] there exist an interval $[a,b]\subset I$, a function $\Phi_1:I\rightarrow\mathbb{R}_+$, $\Phi_1\in L^1(I)$,
	and a constant $c_1\in (0,1]$ satisfying
	\[\begin{array}{rll}
		k_1(t,s)&\leq \Phi_1(s)  & \quad \text{for all } t,s\in I, \\ c_1\,\Phi_1(s)&\leq k_1(t,s) & \quad \text{for all } t\in[a,b], \ s\in I;
	\end{array}\]
	\item[$(H_3)$] there exist a function $\Phi_2:I\rightarrow\mathbb{R}_+$, $\Phi_2\in L^1(I)$, 
	and a constant $c_2\in (0,1]$ such that
	\[c_2\,\Phi_2(s)\leq k_2(t,s)\leq \Phi_2(s)  \quad \text{for all } t,s\in I;\]
	\item[$(H_4)$] the functions $f,g:I\times\mathbb{R}_+^2\rightarrow\mathbb{R}_+$ are continuous.
\end{enumerate}	

In order to apply the theory in Section \ref{sec2}, let us consider the Banach space of continuous functions $X=Y=\mathcal{C}(I)$ endowed with the usual maximum norm $\left\|w\right\|_{\infty} =\max_{t\in I}\left|w(t)\right|$ and the cones 
\begin{align*} 
K_1&=\left\{w\in \mathcal{C}(I)\,:\,w(t)\geq 0 \text{ for all } t\in I \text{ and } \min_{t\in[a,b]}w(t)\geq c_1\left\|w\right\|_{\infty} \right\}, \\ 
K_2&=\left\{w\in \mathcal{C}(I)\,:\, \min_{t\in[0,1]}w(t)\geq c_2\left\|w\right\|_{\infty} \right\}.
\end{align*}
Under assumptions $(H_1)$--$(H_4)$, it can be proven by means of standard arguments (see, for instance, \cite{fig_tojo,Inf}) that $T:=(T_1,T_2)$ maps the cone $K:=K_1\times K_2$ into itself and it is completely continuous; 
in fact, 
if we take $t\in[a,b]\subset [0,1]$, let $(\tilde u(t),\tilde v(t))=T((u,v)){(t)}$, we have
\begin{multline*}	
\tilde u(t)=\int_0^1k_1(t,s)f(s,u(s),v(s))ds\underset{H_2}{\ge}c_1\int_0^1{\Phi_1}(s)f(s,u(s),v(s))ds\\
\underset{H_2}{\ge}c_1\int_0^1 k_1(t',s)f(s,u(s),v(s))ds=c_1 \tilde u(t'), \quad \forall t\in[a,b],\,\,\forall t'\in[0,1].
\end{multline*}
A similar argument holds for the second component.
The compactness of the operator $T$ follows by a direct application of the Arzel\`{a}-Ascoli theorem.

For $\rho>0$, we shall make use of the following open set
\[V_{\rho}=\left\{w\in K_1\,:\,\min_{t\in[a,b]}w(t)<\rho \right\}. \]
Observe that $(K_1)_{\rho}\subset V_{\rho}\subset (K_1)_{\rho/c_1}$, where $(K_1)_{\rho}:=\{w\in K_1\,:\,\left\|w\right\|_{\infty}<\rho \}$. The sets of type $V_\rho$ were introduced by Lan in \cite{Lan} and later employed by several authors, see \cite{Inf} and the references therein.

We are in a position to establish an existence result for the system of Hammerstein type equations \eqref{eq_Ham} as a consequence of Theorem \ref{th_KrasSchauder}.

\begin{theorem}\label{th_Ham}
	Assume that conditions $(H_1)-(H_4)$ are fulfilled. Moreover, suppose that there exist positive numbers $\rho_1,\rho_2>0$, with $\rho_1/c_1<\rho_2$ (resp., $\rho_2<\rho_1$), and $0<\alpha<\beta$ such that
	\begin{enumerate}
		\item[$(H_5)$] there exists a continuous function $\underline{f}:I\rightarrow\mathbb{R}_{+}$ such that
		\[\underline{f}(t)\leq f(t,u,v) \text{ on } [a,b]\times [\rho_1,\rho_1/c_1]\times [\alpha,\beta] \]
		and 
		\[\min_{t\in[a,b]}\int_{a}^{b}k_1(t,s)\underline{f}(s)\,ds\geq \rho_1; \]
		\item[$(H_6)$] there exists a continuous function $\overline{f}:I\rightarrow\mathbb{R}_{+}$ such that
		\[f(t,u,v)\leq \overline{f}(t) \text{ on } [0,1]\times[0,\rho_2]\times[\alpha,\beta] \]
		and
		\[\max_{t\in [0,1]}\int_{0}^{1}k_1(t,s)\overline{f}(s)\,ds\leq \rho_2; \]
		\item[$(H_7)$] there exists a continuous function $\underline{g}:I\rightarrow\mathbb{R}_+$ such that
		\[\underline{g}(t)\leq g(t,u,v) \text{ on } [a,b]\times[\rho_1,\rho_2]\times[\alpha,\beta] \quad (\text{resp., on } [a,b]\times[c_1\,\rho_2,\rho_1/c_1]\times[\alpha,\beta]) \]
		and
		\[\min_{t\in[0,1]}\int_{a}^{b}k_2(t,s)\underline{g}(s)\,ds\geq \alpha;  \]
		\item[$(H_8)$] there exists a continuous function $\overline{g}:I\rightarrow\mathbb{R}_+$ such that
		\[ g(t,u,v)\leq \overline{g}(t) \text{ on } [0,1]\times[0,\rho_2]\times[\alpha,\beta] \quad (\text{resp., on } [0,1]\times[0,\rho_1/c_1]\times[\alpha,\beta]) \]
		and
		\[\max_{t\in[0,1]}\int_{0}^{1}k_2(t,s)\overline{g}(s)\,ds\leq \beta.  \]
	\end{enumerate} 	
	Then the system \eqref{eq_Ham} has at least one positive solution $(u,v)\in K$ such that $\rho_1\leq\displaystyle\min_{t\in[a,b]}u(t)$, $\left\|u\right\|_{\infty}\leq\rho_2$ (resp., $\rho_2\leq\left\|u\right\|_{\infty}$, $\displaystyle\min_{t\in[a,b]}u(t)\leq\rho_1$) and $\alpha\leq v(t)\leq\beta$ for all $t\in I$.
\end{theorem}

\noindent
{\bf Proof.} Suppose that $\rho_1/c_1<\rho_2$ and consider the integral operator $T=(T_1,T_2):\left((\overline{K}_1)_{\rho_2}\setminus V_{\rho_1}\right)\times[\hat{\alpha},\hat{\beta}]\rightarrow K$ defined as above, where $\hat{\alpha}$ and $\hat{\beta}$ denote the constant functions $\hat{\alpha}(t)=\alpha$ and $\hat{\beta}(t)=\beta$ for all $t\in I$. In the sequel, with abuse of notation, $\hat{\alpha}$ and $\hat{\beta}$ will be simply denoted as $\alpha$ and $\beta$, respectively. Let us apply Corollary \ref{hybrid_FPT} in order to obtain a positive solution for the system \eqref{eq_Ham}.

First, let us check that for every $(u,v)\in \left((\overline{K}_1)_{\rho_2}\setminus V_{\rho_1}\right)\times[\alpha,\beta]$ the following conditions are satisfied:
\begin{enumerate}
	\item[$1)$] $T_1(u,v)+\mu\,{\hat {1}}\neq u$ if $u\in \partial\,V_{\rho_1}$ and $\mu>0$ (where $\hat {1}$ denotes the constant function equal to one, {note that $\hat {1}\in K_1$}); 
	\item[$2)$] $T_1(u,v)\neq \lambda\, u$ if $\left\|u\right\|_{\infty}=\rho_2$ and $\lambda> 1$.
\end{enumerate}

To verify that $1)$ holds, assume on the contrary that there exist $(u,v)$ belonging to $\left((\overline{K}_1)_{\rho_2}\setminus V_{\rho_1}\right)\times[\alpha,\beta]$ with $\min_{t\in[a,b]}u(t)=\rho_1$ and $\mu>0$ such that $T_1(u,v)+\mu\,\hat {1}= u$. Then we have 
\[u(t)=\displaystyle\int_{0}^{1} k_1(t,s)f(s,u(s),v(s))\,ds+\mu. \]
Since $u\in K_1$ with $\min_{t\in[a,b]}u(t)=\rho_1$, it follows from the definition of the cone $K_1$ that $\rho_1\leq u(t)\leq \rho_1/c_1$ for all $t\in[a,b]$. Hence, for any $t\in[a,b]$, we deduce from hypothesis $(H_5)$ that
\[u(t)> \int_{a}^{b} k_1(t,s)f(s,u(s),v(s))\,ds\geq \int_{a}^{b} k_1(t,s)\underline{f}(s)\,ds\geq \rho_1, \]
a contradiction.

Now, let us show that $\left\|T_1(u,v)\right\|_{\infty}\leq\rho_2$ for all $(u,v)\in K_1\times[\alpha,\beta]$ with $\left\|u\right\|_{\infty}=\rho_2$, which clearly implies property $2)$. Note that, for such $(u,v)$, we have $0\leq u(t)\leq \rho_2$ and $\alpha\leq v(t)\leq \beta$ for all $t\in I$. Thus, we get, for $t\in I$,
\[T_1(u,v)(t)=\displaystyle\int_{0}^{1} k_1(t,s)f(s,u(s),v(s))\,ds\leq\int_{0}^{1} k_1(t,s)\overline{f}(s)\,ds. \]
By condition $(H_6)$, we deduce $\left\|T_1(u,v)\right\|_{\infty}\leq \max_{t\in[0,1]}\int_{0}^{1} k_1(t,s)\overline{f}(s)\,ds\leq\rho_2$.

Finally, it remains to show that $T_2\left(\left((\overline{K}_1)_{\rho_2}\setminus V_{\rho_1}\right)\times[\alpha,\beta] \right)\subset[\alpha,\beta]$. Take $(u,v)\in \left((\overline{K}_1)_{\rho_2}\setminus V_{\rho_1}\right)\times[\alpha,\beta]$ and let us check that $\alpha\leq T_2(u,v)(t)\leq \beta$ for all $t\in I$. On the one hand, being $k_2$ and $g$ non-negative functions, we have that for $t\in[0,1]$,
\[T_2(u,v)(t)=\int_{0}^{1} k_2(t,s)g(s,u(s),v(s))\,ds\geq \int_{a}^{b} k_2(t,s)g(s,u(s),v(s))\,ds \]
and thus, since $\rho_1\leq u(t)\leq \rho_2$ for all $t\in[a,b]$, condition $(H_7)$ implies that
\[T_2(u,v)(t)\geq \min_{t\in[0,1]} \int_{a}^{b} k_2(t,s)\underline{g}(s)\,ds\geq \alpha. \]
On the other hand, by condition $(H_8)$, we have that for $t\in[0,1]$,
\[T_2(u,v)(t)\leq \max_{t\in[0,1]}\int_{0}^{1} k_2(t,s)g(s,u(s),v(s))\,ds\leq \max_{t\in[0,1]}\int_{0}^{1} k_2(t,s)\overline{g}(s)\,ds\leq\beta. \]

In conclusion, Corollary \ref{hybrid_FPT} ensures that the operator $T$ has at least one fixed point localized in the set $\left((\overline{K}_1)_{\rho_2}\setminus {V}_{\rho_1}\right)\times[\alpha,\beta]$. Finally, note that if $\rho_2<\rho_1$, the proof is analogous, but now considering the operator $T$ defined in the set $\left( \overline{V}_{\rho_1}\setminus (K_1)_{\rho_2}\right)\times[\alpha,\beta]$.
\qed

As an illustrative example, we deal with the existence of positive solutions for the following system of second-order equations
\begin{equation}\label{syst_2or}
	\begin{array}{r}
		u''(t)+f(t,u(t),v(t))=0, \quad t\in[0,1], \\[0.3cm]
		v''(t)+g(t,u(t),v(t))=0, \quad t\in[0,1],
	\end{array}
\end{equation}
coupled with the following boundary conditions
\begin{equation}\label{BCs}
	u(0)=u(1)=0=v'(0)=v(1)+v'(1).
\end{equation}	

The system \eqref{syst_2or}--\eqref{BCs} can be studied by means of a system of integral equations of the form \eqref{eq_Ham}, where the kernels $k_1$ and $k_2$ are the corresponding Green's functions, which are given by
\[k_1(t,s)=\left\{\begin{array}{ll} {(1-t)s}, & \quad s\leq t, \\ {t(1-s)}, & \quad s>t, \end{array} \right. \]
and 
\[k_2(t,s)=\left\{\begin{array}{ll} 2-t, & \quad s\leq t, \\ 2-s, & \quad s>t. \end{array} \right. \]
It is well--known (see, for example, \cite{Inf}) that the Green's function $k_1$ satisfies condition $(H_2)$ if we take 
\[\Phi_1(s)=s(1-s), \ s\in I, \quad [a,b]=[1/4,3/4], \quad c_1=1/4.  \]
On the other hand, the Green's function $k_2$ was studied in \cite{Webb}  (see also \cite{CCI14}) and it was shown that it satisfies condition $(H_3)$ with \[\Phi_2(s)=2-s, \ s\in I,\quad c_2=1/2.\] 

\begin{example}\label{numex}
	Consider the system \eqref{syst_2or}--\eqref{BCs} with the continuous nonlinearities
	\[f(t,u,v)=t\,u^2(1+\sin^2 v), \quad g(t,u,v)=t(2+\sin u)(6+\cos v). \]
	Take $\rho_1=64$, $\rho_2=4$, $\alpha=1$ and {$\beta=14$. Then} it is easy to check that hypotheses $(H_5)-(H_8)$ hold with $\underline{f}(t)=1024$,  $\overline{f}(t)=32$, $\underline{g}(t)=5\,t$ and $\overline{g}(t)=21\,t$, $t\in I$.
	
	Therefore, Theorem \ref{th_Ham} guarantees that the system \eqref{syst_2or}--\eqref{BCs} has a positive solution. We emphasize that the theory developed in \cite{CCI14} cannot be applied to this example since the function $g(t,\cdot,\cdot)$ is not monotone in any rectangle of the form $[0,B_1]\times[0,B_2]$.

We run a numerical approximation for the solution $(u,v)$ of the system, with the aid of the  MATLAB Function BVP4c; the result is illustrated in Figure~\ref{fig:test}. We started with the initial guess
$$u_0(t)\approx 50.667 t -99.333 t^2 + 85.333 t^3 -42.667 t^4, \quad v_0 (t)\approx 4 -0.2 t -1.6 t^2,$$
based on the properties of the desired solution.
The numerically obtained solution is consistent with the theoretical result, in the sense that
the two components of the solution have the following properties predicted by theory:
\begin{itemize}
\item $u\ge 0$, $\displaystyle\min_{t\in[1/4,3/4]}u(t)\ge 1/4\left\|u\right\|_{\infty}$, $4\leq\left\|u\right\|_{\infty}$, $\displaystyle\min_{t\in[1/4,3/4]}u(t)\leq 64$, 	
\item $1\le v(t)\le 14$,
\end{itemize}	
as can be seen from {Figure~\ref{fig:test}}.
	
\begin{figure}[h]
\centering
\begin{subfigure}{.5\textwidth}
  \centering
  \includegraphics[width=.7\linewidth]{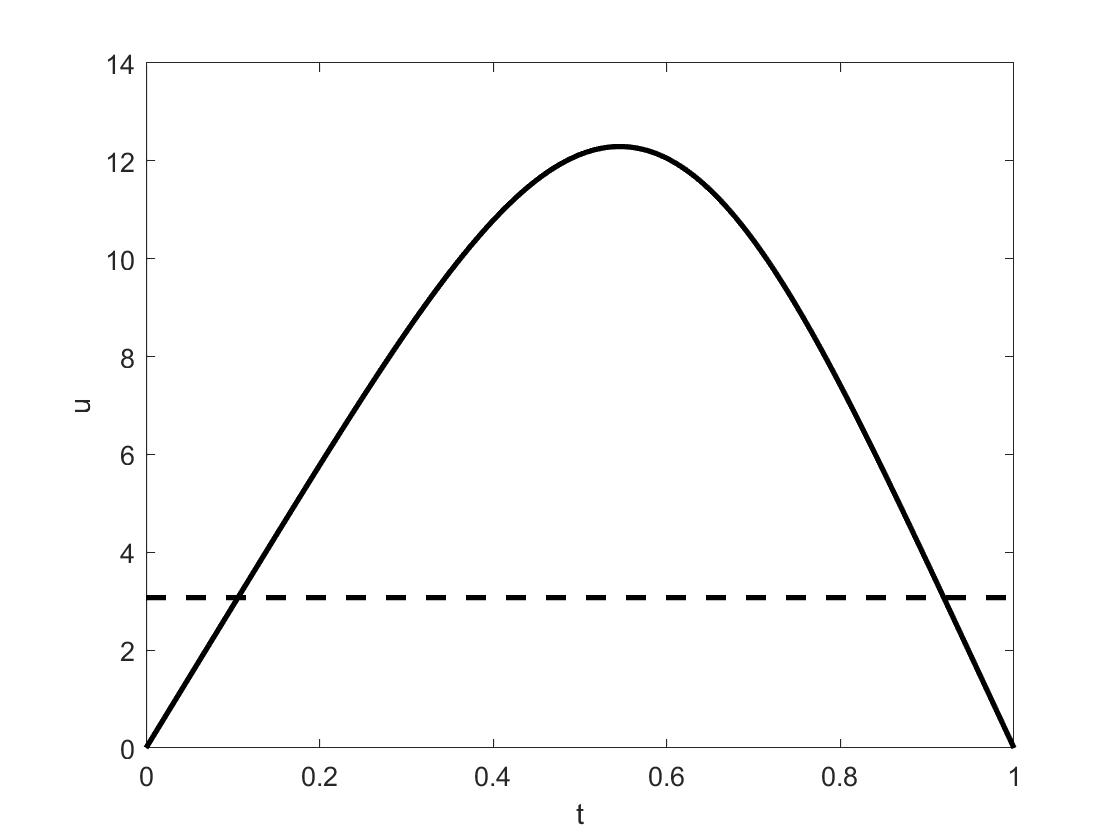}
  \caption{The component $u$; the dashed line represents $1/4\left\|u\right\|_{\infty}$.}
\end{subfigure}%
\begin{subfigure}{.5\textwidth}
  \centering
  \includegraphics[width=.7\linewidth]{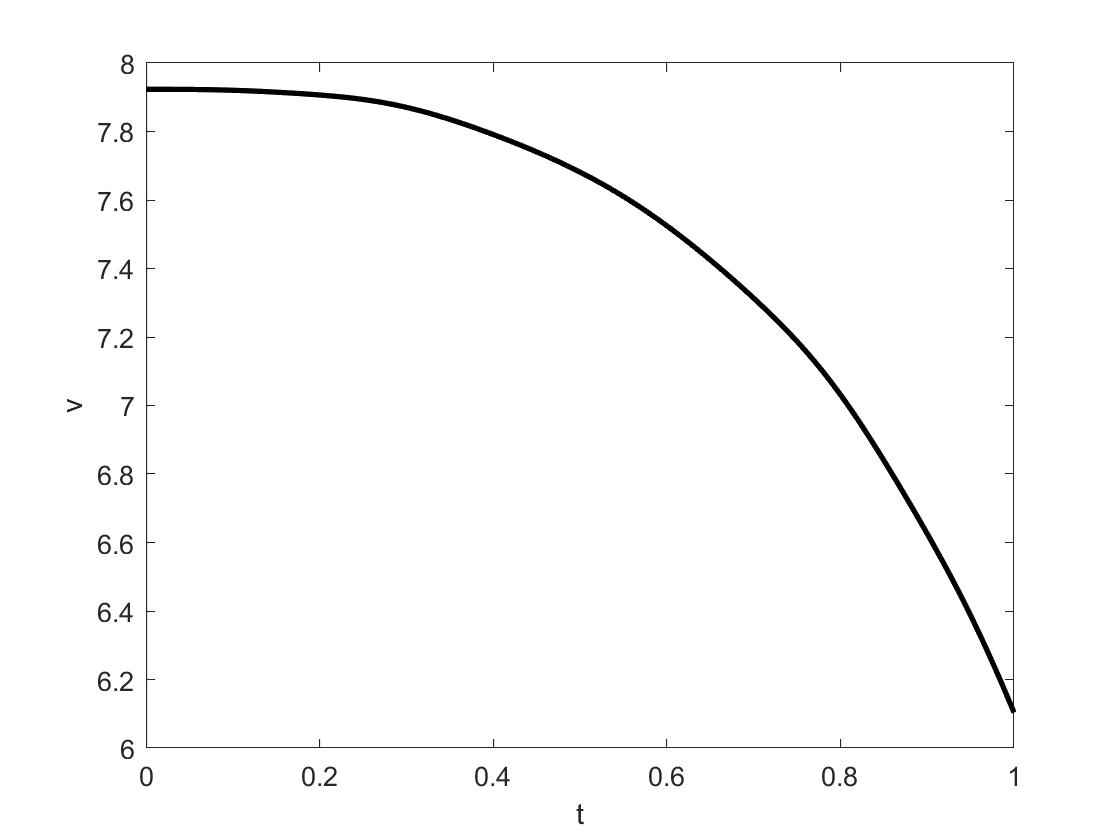}
  \caption{The component $v$.}
\end{subfigure}
\caption{A numerical solution for Example~\ref{numex}.}
\label{fig:test}
\end{figure}
\end{example}

It is worth mentioning that Theorem \ref{th_KrasSchauder} allows us to obtain other kind of localization of the solutions of the system \eqref{eq_Ham} depending on the closed convex set $D$ we choose. Of course, this will be related to the behavior of the nonlinearity $g$. In order to illustrate this fact, let us choose $D$ as a ball in the following existence result concerning system \eqref{eq_Ham}.

\begin{theorem}\label{th_Ham2}
	Under hypotheses $(H_1)-(H_3)$, let us assume that the functions $f:I\times\mathbb{R}_{+}\times\mathbb{R}\rightarrow\mathbb{R}_{+}$ and $g:I\times\mathbb{R}_{+}\times\mathbb{R}\rightarrow\mathbb{R}$ are continuous and there exist $\rho_1,\rho_2>0$, with $\rho_1/c_1<\rho_2$ (resp., $\rho_2<\rho_1$), and $R_2>0$ such that the following conditions are satisfied:
	\begin{enumerate}
		\item[$(H_5^*)$] there exists a continuous function ${f}_*:I\rightarrow\mathbb{R}_{+}$ such that
		\[{f}_*(t)\leq f(t,u,v) \text{ on } [a,b]\times [\rho_1,\rho_1/c_1]\times [-R_2,R_2] \]
		and 
		\[\min_{t\in[a,b]}\int_{a}^{b}k_1(t,s){f}_*(s)\,ds\geq \rho_1; \]
		\item[$(H_6^*)$] there exists a continuous function ${f}^*:I\rightarrow\mathbb{R}_{+}$ such that
		\[f(t,u,v)\leq {f}^*(t) \text{ on } [0,1]\times[0,\rho_2]\times[-R_2,R_2] \]
		and
		\[\max_{t\in [0,1]}\int_{0}^{1}k_1(t,s){f}^*(s)\,ds\leq \rho_2; \]
		\item[$(H_7^*)$] there exists a continuous function ${g}^*:I\rightarrow\mathbb{R}_{+}$ such that
		\[\left|g(t,u,v) \right|\leq g^*(t) \text{ on } [0,1]\times[0,R_1]\times[-R_2,R_2] \quad (\text{where } R_1:=\max\{\rho_1/c_1,\rho_2 \}) \]
		and 
		\[\max_{t\in [0,1]}\int_{0}^{1}k_2(t,s){g}^*(s)\,ds\leq R_2. \]
	\end{enumerate}	
	Then the system \eqref{eq_Ham} has at least one non-trivial solution $(u,v)\in K_1\times X$ such that $\rho_1\leq\displaystyle\min_{t\in[a,b]}u(t)$, $\left\|u\right\|_{\infty}\leq\rho_2$ (resp., $\rho_2\leq\left\|u\right\|_{\infty}$, $\displaystyle\min_{t\in[a,b]}u(t)\leq\rho_1$) and $\left\|v\right\|_{\infty}\leq R_2$.
\end{theorem}

\noindent
{\bf Proof.} It follows as an application of Theorem \ref{th_KrasSchauder} to the operator $T=(T_1,T_2)$ defined above and the choice of the set $D$ as the closed ball of radius $R_2$, that is, \[D=\{w\in\mathcal{C}(I)\,:\,\left\|w\right\|_{\infty}\leq R_2 \}.\] The details are similar to those in the proof of Theorem \ref{th_Ham}.
\qed

\begin{example}\label{ex2}
	Consider the system \eqref{syst_2or}--\eqref{BCs} with the nonlinearities 
	\[f(t,u,v)=t\,u^2(1+\sin^2 v), \quad g(t,u,v)=t\,e^{v^2-2}\sin u. \]
	A simple computation shows that $\max_{t\in I}\int_{0}^{1}k_2(t,s)\,ds=3/2$ in this case, so if we take $\rho_1=64$, $\rho_2=4$ and $R_2=1$, then it is easy to check that hypotheses $(H_5^*)-(H_7^*)$ are satisfied.
	
	Therefore, Theorem \ref{th_Ham2} guarantees that the system \eqref{syst_2or}--\eqref{BCs} has at least one solution $(u,v)$ such that $u$ is positive (nontrivial) and $\left|v(t)\right|\leq 1$ for all $t\in I$. Since $g$ is a sign-changing nonlinearity, $v$ is not necessarily positive, but clearly it cannot be the identically zero function.  
	
Also in this case, we found a numerical approximation for one solution $(u,v)$ of the system, with the aid of the same MATLAB Functions; the results are reported in Figure~\ref{fig:test2}. As the initial guess we used the same $u_0$ as in the previous example, and as $v_0$ we chose a function proportional by a factor $10^3$ to the solution
of second differential equation with $g=g(t,u_0,v)$. After that, we iterated the solution starting from the previous one, up to when a relative tolerance not greater than $10^{-7}$ was reached.
The obtained solution is consistent with the theoretical result.
	
\begin{figure}[h]
\centering
\begin{subfigure}{.5\textwidth}
  \centering
  \includegraphics[width=.7\linewidth]{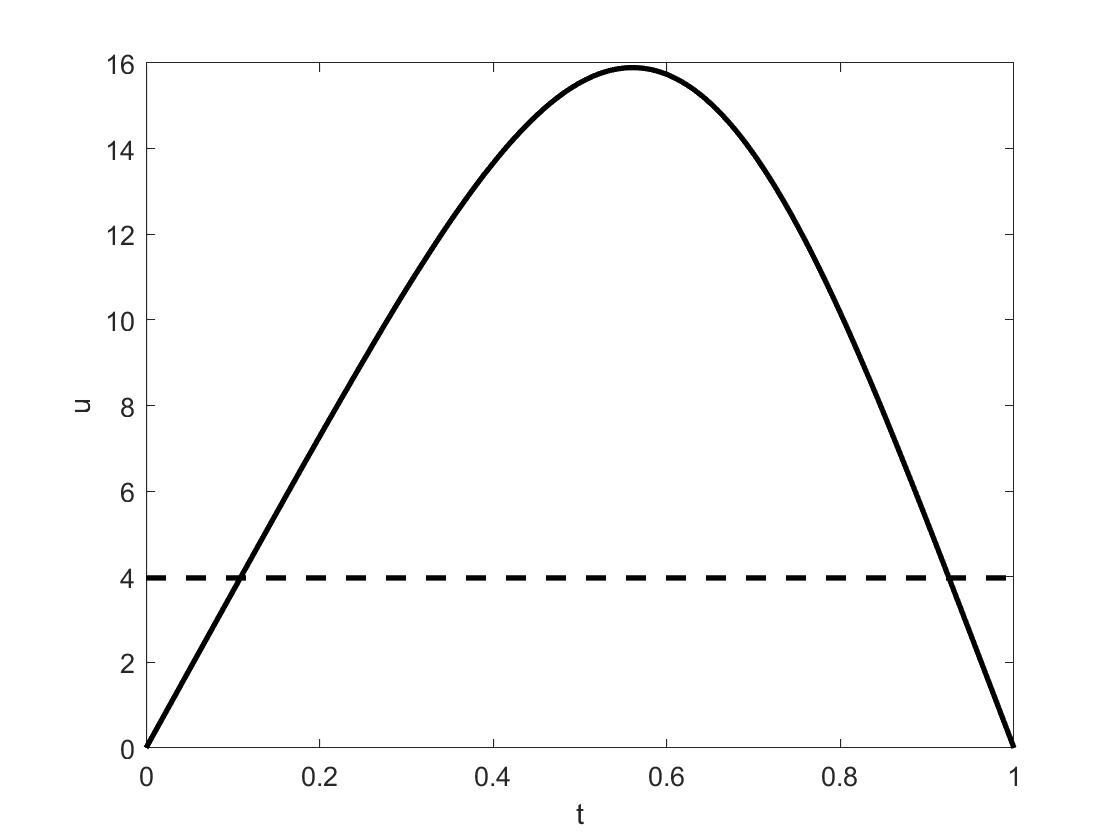}
  \caption{The component $u$; the dashed line represents $1/4\left\|u\right\|_{\infty}$.}
\end{subfigure}%
\begin{subfigure}{.5\textwidth}
  \centering
  \includegraphics[width=.7\linewidth]{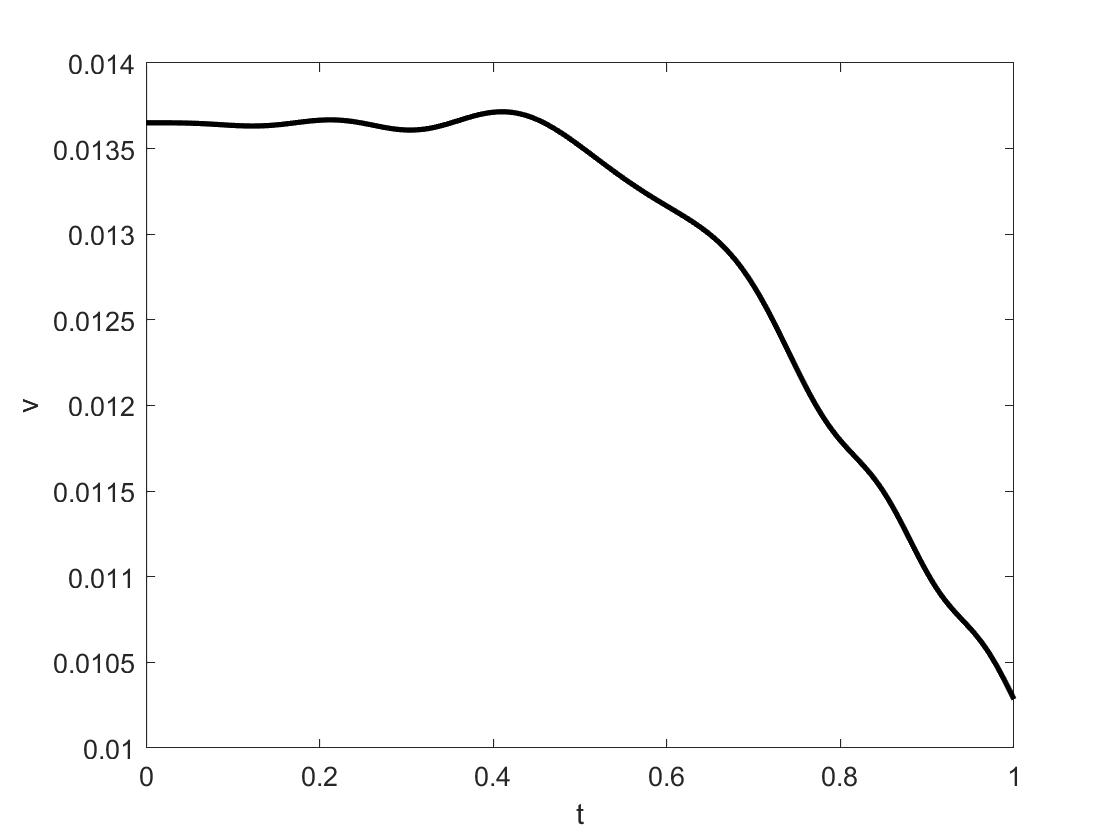}
  \caption{The component $v$.}
\end{subfigure}
\caption{A numerical solution for Example~\ref{ex2}.}
\label{fig:test2}
\end{figure}

\end{example}

\section*{Acknowledgements}
G.~Infante is a member of the Gruppo Nazionale per l'Analisi Matematica, la Probabilit\`a e le loro Applicazioni (GNAMPA) of the Istituto Nazionale di Alta Matematica (INdAM) and of the UMI Group TAA  ``Approximation Theory and Applications''. G.~Mascali is a member of the Gruppo Nazionale per la Fisica Matematica (GNFM) of INdAM. G.~Infante and G.~Mascali are supported by the project POS-CAL.HUB.RIA. 
J. Rodr\'iguez--L\'opez has been partially supported by the VIS Program of the University of Calabria, and by Ministerio de Ciencia y Tecnología (Spain), AIE and Feder, grant PID2020-113275GB-I00. 
This paper was written during a visit of J. Rodr\'iguez-L\'opez to the Dipartimento di Matematica e Informatica of the Universit\`a della Calabria. J. Rodr\'iguez-L\'opez is grateful to the people of the aforementioned Dipartimento for their kind and warm hospitality.

\end{document}